\documentclass{amsart}
\usepackage{graphicx}
\usepackage{latexsym,amsmath,amssymb}
\theoremstyle{definition}
\theoremstyle{remark}

\numberwithin{equation}{section}
\begin{document}

\title[Weighted conditional type operators ]
{ The spectra of weighted conditional type operators }

\author{\sc\bf Y. Estaremi  }
\address{\sc Y. Estaremi  }
\email{yestaremi@pnu.ac.ir}

\address{Department of mathematics, Payame Noor university , p. o. box: 19395-3697, Tehran,
Iran}

\thanks{}

\thanks{}



\date{}

\dedicatory{}

\commby{}

%%% ----------------------------------------------------------------------
\begin{abstract}
In this paper, the spectrum, residual spectrum, point spectrum and
spectral radius of weighted conditional type operators are
computed. As an application, we give an equivalent condition for
weighted conditional type operators to be quasinilpotent. Also,
several examples are provided to illustrate concrete application
of the main results of the paper.

\noindent {}
\end{abstract}

\maketitle

\section{ \sc\bf Introduction }
This paper is about an important operator in statistics and
analysis, that is called conditional expectation. Theory of
conditional type operators is one of important arguments in the
connection of operator theory and measure theory. By the
projection theorem, the conditional expectation $E(X)$ is the best
mean square predictor of $X$ in range $E$, also in analysis it is
proved that lots of operators are of the form $E$ and of the form
of combinations of $E$ and multiplications operators. I think our
information in this paper and references \cite{ej1, ej2, e} will
be useful for experts in statistics. Conditional expectations have
been studied in an operator theoretic setting, by, for example,in
\cite{mo}, S.-T. C. Moy characterized all operators on $L^p$ of
the form $f\rightarrow E(fg)$ for $g$ in $L^q$ with $E(|g|)$
bounded, P.G. Dodds, C.B. Huijsmans and B. De Pagter, \cite{dhd},
extended these characterizations to the setting of function ideals
and vector lattices and J. Herron presented some assertions about
the operator $EM_u$ on $L^p$ spaces in \cite{her}. Also, some
results about multiplication conditional expectation operators can
be found in \cite{g, lam}. In \cite{e,ej1, ej2} we investigated
some classic properties of multiplication conditional expectation
operators $M_wEM_u$ on $L^p$ spaces.\\

\section{ \sc\bf Preliminaries}
Let $(X,\Sigma,\mu)$ be a complete $\sigma$-finite measure space.
For any sub-$\sigma$-finite algebra $\mathcal{A}\subseteq
 \Sigma$, the $L^2$-space
$L^2(X,\mathcal{A},\mu_{\mid_{\mathcal{A}}})$ is abbreviated  by
$L^2(\mathcal{A})$, and its norm is denoted by $\|.\|_2$. All
comparisons between two functions or two sets are to be
interpreted as holding up to a $\mu$-null set. The support of a
measurable function $f$ is defined as $S(f)=\{x\in X; f(x)\neq
0\}$. We denote the vector space of all equivalence classes of
almost everywhere finite valued measurable functions on $X$ by
$L^0(\Sigma)$.

\vspace*{0.3cm} For a sub-$\sigma$-finite algebra
$\mathcal{A}\subseteq\Sigma$, the conditional expectation operator
associated with $\mathcal{A}$ is the mapping $f\rightarrow
E^{\mathcal{A}}f$, defined for all non-negative measurable
function $f$ as well as for all $f\in L^2(\Sigma)$, where
$E^{\mathcal{A}}f$, by the Radon-Nikodym theorem, is the unique
$\mathcal{A}$-measurable function satisfying
$$\int_{A}fd\mu=\int_{A}E^{\mathcal{A}}fd\mu, \ \ \ \forall A\in \mathcal{A} .$$
As an operator on $L^{2}({\Sigma})$, $E^{\mathcal{A}}$ is
idempotent and $E^{\mathcal{A}}(L^2(\Sigma))=L^2(\mathcal{A})$.
This operator will play a major role in our work. Let $f\in
L^0(\Sigma)$, then $f$ is said to be conditionable with respect to
$E$ if $f\in\mathcal{D}(E):=\{g\in L^0(\Sigma): E(|g|)\in
L^0(\mathcal{A})\}$. Throughout this paper we take $u$ and $w$ in
$\mathcal{D}(E)$. If there is no possibility of confusion, we
write $E(f)$ in place of $E^{\mathcal{A}}(f)$.  A detailed
discussion about this operator may be found in \cite{rao}.\\

Let $\mathcal{H}$ be the infinite dimensional complex Hilbert
space and let $\mathcal{L(H)}$ be the algebra of all bounded
operators on $\mathcal{H}$. Every operator $T$ on a Hilbert space
$\mathcal{H}$ can be decomposed into $T = U|T|$ with a partial
isometry $U$, where $|T| = (T^*T)^{\frac{1}{2}}$ . $U$ is
determined uniquely by the kernel condition $\mathcal{N}(U) =
\mathcal{N}(|T|)$. Then this decomposition is called the polar
decomposition. The Aluthge transformation $\widehat{T}$ of the
operator $T$ is defined by
$\widehat{T}=|T|^{\frac{1}{2}}U|T|^{\frac{1}{2}}$. \\
In this paper we will be concerned with characterizing weighted
conditional expectation type operators on $L^2(\Sigma)$, computing
the spectrum, residual spectrum, point spectrum, spectral radius.

 \section{ \sc\bf The spectra of weighted conditional type operators}

In the first we reminisce some properties of weighted conditional
type operators that we proved in
\cite{ej1}. Also, we give some example of conditional expectation operator.\\

Let $T=M_wEM_u$ be a bounded operator on $L^{2}(\Sigma)$ and let
$p\in (0,\infty)$. Then
$$(T^{\ast}T)^{p}=M_{\bar{u}(E(|u|^{2}))^{p-1}\chi_{S}(E(|w|^{2}))^{p}}EM_{u}$$
and
$$(TT^{\ast})^{p}=M_{w(E(|w|^{2}))^{p-1}\chi_{G}(E(|u|^{2}))^{p}}EM_{\bar{w}},$$
where $S=S(E(|u|^2))$ and $G=S(E(|w|^2))$. And the unique polar
decomposition of bounded operator $T=M_wEM_u$ is $U|T|$, where

$$|T|(f)=\left(\frac{E(|w|^{2})}{E(|u|^{2})}\right)^{\frac{1}{2}}\chi_{S}\bar{u}E(uf)$$
and
 $$U(f)=\left(\frac{\chi_{S\cap
 G}}{E(|w|^{2})E(|u|^{2})}\right)^{\frac{1}{2}}wE(uf),$$
for all $f\in L^{2}(\Sigma)$. Also, the Aluthge transformation of
$T=M_wEM_u$ is
$$\widehat{T}(f)=\frac{\chi_{S}E(uw)}{E(|u|^{2})}\bar{u}E(uf), \ \ \ \ \ \  \ \ \ \  \ \  \ f\in
L^{2}(\Sigma).$$\\

From now on, we shall denote by $\sigma(T)$, $\sigma_{p}(T)$,
$\sigma_{jp}(T)$, $\sigma_r(T)$, $r(T)$ the spectrum of $T$, the
point spectrum of $T$, the joint point spectrum of $T$, the
residual spectrum, the spectral radius of $T$, respectively. The
spectrum of an operator $T$ is the set
$$\sigma(T)=\{\lambda\in \mathbb{C}:T-\lambda I \ \ \  {\text{is \  not \
invertible}}\}.$$ A complex number $\lambda\in \mathbb{C}$ is said
to be in the point spectrum $\sigma_{p}(T)$ of the operator $T$,
if there is a unit vector $x$ satisfying $(T-\lambda)x=0$. If in
addition, $(T^{\ast}-\bar{\lambda})x=0$, then $\lambda$ is said to
be in the joint point spectrum $\sigma_{jp}(T)$ of $T$. The
residual spectrum of $T$ is equal to $\mathbb{C} \setminus
(\sigma_p(EM_u)\cup \rho(EM_u))$. Also, the spectral radius of $T$ is defined by $r(T)=\sup\{|\lambda|: \lambda\in \sigma(T)\}$.\\
If $A, B\in \mathcal{B}(\mathcal{H})$, then it is well known that
$\sigma(AB)\setminus\{0\}=\sigma(BA)\setminus\{0\}$,
$\sigma_{p}(AB)\setminus\{0\}=\sigma_{p}(BA)\setminus\{0\}$,
$\sigma_{jp}(AB)\setminus\{0\}=\sigma_{jp}(BA)\setminus\{0\}$  and
$\sigma_{r}(AB)\setminus\{0\}=\sigma(BA)\setminus\{0\}$.
 J. Herron
showed that if $EM_u:L^2(\Sigma)\rightarrow L^2(\Sigma)$, then
$\sigma(EM_u)=ess \ range(E(u))\cup\{0\}$\cite{her}. By means of
the above mentioned properties of weighted conditional expectation
type operators we have the following theorems.\\

 \vspace*{0.3cm}
{\bf Theorem. 3.1.} Let $EM_u:L^p(\Sigma)\rightarrow
L^p(\mathcal{A})$, for $1\leq p \leq \infty$. Then

$$\rho(EM_{u})=\{\lambda\in \mathbb{C}\mid \exists \epsilon>0 \ \
\ such \ that \ | \lambda- E(u)|\geq \epsilon \ a.e\}$$

$$\sigma(EM_u)=\{\lambda\in \mathbb{C}\mid \nexists \epsilon>0 \ \
\ such \ that \ | \lambda- E(u)|\geq \epsilon \ a.e\}$$

$$\sigma_{p}(EM_u)=\{\lambda\in \mathbb{C}\mid \mu(\{x\in X |
E(u)(x)=\lambda\})>0\}$$\\

$\sigma_r(EM_u)=
\left\{%
\begin{array}{ll}
    \emptyset, & 1\leq p <\infty; \\
    \{\lambda\in \mathbb{C}\mid \nexists \epsilon>0 \ \
\ such \ that \ | \lambda- E(u)|\geq \epsilon \ a.e\}, & p=\infty \\
\end{array}%
\right.$\\

\vspace*{0.3cm} {\bf Proof.} Let $\lambda\in
\mathbb{C}\setminus\{0\}$ such that $|\lambda-E(u)|>\epsilon$,
a.e, for some $\epsilon>0$. If $\lambda f-E(uf)=0$, then $f$ is
$\mathcal{A}$-measurable and so, $(\lambda-E(u))f=0$. This implies
that $f=0$ a.e, therefore $\lambda I-EM_u$ is injective. Hence
$$\lambda \in \mathbb{C}\setminus (\sigma_p(EM_u)\cup \{0\})=\rho(EM_u).$$\\

Conversely, let $\lambda\in \mathbb{C}\setminus
(\sigma_p(EM_u)\cup\{0\})$, then $\lambda I-EM_u$ has an inverse
operator. Define linear transformation $S$ on
$\mathcal{R}(EM_u-\lambda I)\subseteq L^p(\Sigma)$ as follows

$$Sf=\frac{E(uf)-f(E(u)-\lambda)}{\lambda(E(u)-\lambda)}, \ \ \ f\in \mathcal{R}(EM_u-\lambda I).$$

If  there exists $\epsilon>0$ such that $|\lambda -E(u)|>\epsilon$
a. e, then
$\|\frac{1}{E(u)-\lambda}\|_{\infty}\leq\frac{1}{\varepsilon}$. So
$$\|Sf\|_p\leq \|\frac{E(uf)}{\lambda(E(u)-\lambda)}\|_p+\|\frac{f}{\lambda}\|_p$$

$$\leq(\frac{\|EM_u\|}{\lambda \varepsilon}+\frac{1}{\lambda})\|f\|_p.$$
Thus $S$ is bounded $\mathcal{R}(EM_u-\lambda I)$. If $S$ is
bounded on $\mathcal{R}(EM_u-\lambda I)$, then for $f\in
L^p(\mathcal{A})\cap \mathcal{R}(EM_u-\lambda I)$, $Sf=\alpha
f=M_{\alpha}f$, where $\alpha=\frac{1}{Eu-\lambda}$. Thus
multiplication operator $M_{\alpha}$ is bounded on
$L^p(\mathcal{A})\cap \mathcal{R}(EM_u-\lambda I)$. This implies
that $\alpha\in L^{\infty}(\mathcal{A})$ and so there exist some
$\varepsilon>0$ such that
$E(u)-\lambda=\frac{1}{\alpha}\geq\varepsilon$ a.e.
 Also, we have

$$S\circ (EM_u-\lambda I)=I,$$
Thus $(T-\lambda I)^{-1}=S$ and so $S$ have to be bounded. Hence
there exist some $\varepsilon>0$ such that
$E(u)-\lambda=\frac{1}{\alpha}\geq\varepsilon$ a.e.
\\

Let $A_{\lambda}=\{x\in X:E(u)(x)=\lambda\}$, for $\lambda\in
\mathbb{C}$. Suppose that $\mu(A_{\lambda})>0$. Since
$\mathcal{A}$ is $\sigma$-finite, there exists an
$\mathcal{A}$-measurable subset $B$ of $A_{\lambda}$ such that
$0<\mu(B)<\infty$, and $f=\chi_{B}\in L^p(\mathcal{A})\subseteq
L^p(\Sigma)$. Now
$$EM_u(f)-\lambda f=E(u)\chi_{B}-\lambda \chi_{B}=0.$$ This
implies that $\lambda\in \sigma_{p}(EM_u)$.\\
If there exists $f\in L^p(\Sigma)$ such that $f\chi_{C}\neq 0$
$\mu$-a.e, for $C\in \Sigma$ of positive measure and
$E(uf)=\lambda f$ for $\lambda \in \mathbb{C}$, which means that
$f$ is $\mathcal{A}$-measurable. Therefore $E(uf)=E(u)f=\lambda f$
and $(E(u)-\lambda)f=0$. This implies that $C\subseteq
A_{\lambda}$ and so $\mu(A_{\lambda})>0$.\\

Residual spectrum: Let $\lambda\in \mathbb{C} \setminus
(\sigma_p(EM_u)\cup \rho(EM_u))$. So $\mu(\{x\in X |
E(u)(x)=\lambda\})=0$, but on the other hand $\mu(\{x\in X |
|E(u)(x)-\lambda|>\epsilon\})>0$ for every $\epsilon>0$. We wish
to determine if the range of $\lambda I-EM_u$, i.e. the domain of
$(\lambda I-EM_u)^{-1}$, is dense. Set, for each $n\in
\mathbb{N}$,
$$E_n =\{x\in X: |E(u)(x)-\lambda|\geq\frac{1}{n}\}$$ The range of
$\lambda I-EM_u$ contains $\{\chi_{E_n}f:f\in L^p(\mathcal{A}),
n\in \mathbb{N}\}$, because, for every $f\in L^p(\mathcal{A})$,
$(\lambda I-EM_u)(\frac{1}{\lambda-E(u)}\chi_{E_n}
f)=\chi_{E_n}f$.  Furthermore, $\chi_{E_n}f$ converges pointwise
almost everywhere to $f$ as $n\rightarrow \infty$. So if $1\leq
p<\infty$, the Lebesgue dominated convergence theorem implies that
$\chi_{E_n}f$ converges in $L^p(X, \mathcal{A}, \mu)$ to $f$ as
$n\rightarrow \infty$. Thus the range of range of $\lambda I-EM_u$
is dense in $L^p(X, \mathcal{A}, \mu)$ if $1\leq p <\infty$. On
the other hand, if $p=\infty$, the constant function $1$ is not in
the closure of the range of $\lambda I-EM_u$ because, for every
$0\neq f\in L^p(X, \mathcal{A}, \mu)$, there is $A\in \mathcal{A}$
such that $\mu(A)>0$ and  $|\lambda -E(u)|\leq \frac{1}{2
\|f\|_{\infty}}$ on $A$ and hence $|1- (\lambda- E(u))f)|\geq
\frac{1}{2}$ on $A$. Thus the proof is completed.\\

 \vspace*{0.3cm} {\bf Theorem 3.2.}
Let
$T=M_wEM_u:L^p(\Sigma)\rightarrow L^p(\Sigma)$, for $1\leq p\leq \infty$. Then\\

(a) $$\sigma(M_wEM_u)\setminus\{0\}=\{\lambda\in \mathbb{C}\mid
\nexists \epsilon>0 \ \ \ such \ that \ | \lambda- E(uw)|\geq
\epsilon \ a.e\}\setminus\{0\} $$

(b) If $S\cap G=X$ and $p=2$, then $$\sigma(M_wEM_u)=\{\lambda\in
\mathbb{C}\mid \nexists \epsilon>0 \ \ \ such \ that \ | \lambda-
E(uw)|\geq \epsilon \ a.e\},$$ where $S=S(E(|u|^2))$ and
$G=S(E(|w|^2))$.

(c)
$$\sigma_{p}(M_wEM_u)\setminus\{0\}=\{\lambda\in\mathbb{C}\setminus\{0\}:\mu(A_{\lambda,w})>0\},$$
where $A_{\lambda,w}=\{x\in X:E(uw)(x)=\lambda\}$.\\

(d) $$\rho(M_wEM_{u})\setminus\{0\}=\{\lambda\in \mathbb{C}\mid
\exists \epsilon>0 \ \ \ such \ that \ | \lambda- E(uw)|\geq
\epsilon \ a.e\}\setminus\{0\}.$$

(e) $\sigma_r(M_wEM_u)\setminus\{0\}=
\left\{%
\begin{array}{ll}
    \emptyset, & 1\leq p <\infty; \\
    \{\lambda\in \mathbb{C}\mid \nexists \epsilon>0 \ \
\ such \ that \ | \lambda- E(uw)|\geq \epsilon \ a.e\}\setminus\{0\}, & p=\infty \\
\end{array}%
\right.$\\

 \vspace*{0.3cm} {\bf Proof.} (a) Since
 $$\sigma(M_wEM_u)\setminus\{0\}=\sigma(EM_uM_w)\setminus\{0\}=\sigma(EM_{uw})\setminus\{0\}=ess
 \ range(E(uw))\setminus\{0\},$$
 then we have
 $$\sigma(M_wEM_u)\setminus \{0\}=ess \
range(E(uw))\setminus\{0\}.$$

(b) We know that $\sigma(EM_{uw})=ess \ range(E(uw))$. So, we have
to prove that $0\notin\sigma(EM_{uw})$ if and only if
$0\notin\sigma(M_wEM_u)$ by (a).\\
Let $0\notin\sigma(EM_{uw})$. Then $EM_{uw}$ is surjective and so
$\mathcal{A}=\Sigma$. Thus $E=I$. So $0\notin\sigma(M_wEM_u)$.\\
Conversely, we know that the polar decomposition of
$M_wEM_u=U|M_wEM_u|$ is as follow,
$$|M_wEM_u|(f)=\left(\frac{E(|w|^{2})}{E(|u|^{2})}\right)^{\frac{1}{2}}\chi_{S}\bar{u}E(uf)$$
and
 $$U(f)=\left(\frac{\chi_{S\cap
 G}}{E(|w|^{2})E(|u|^{2})}\right)^{\frac{1}{2}}wE(uf),$$
for all $f\in L^{2}(\Sigma)$.\\

If $0\notin\sigma(M_wEM_u)$, then $|M_wEM_u|$ is invertible and
$U$ is unitary. Therefore $U^{\ast}U=UU^{\ast}=I$. The equation
$UU^{\ast}=I$ implies that $w\in L^{0}(\mathcal{A}_{S\cap G})$,
where $\mathcal{A}_{S\cap G}=\{A\cap S\cap G: \ A\in
\mathcal{A}\}$. Since $S\cap G=X$, then $w\in L^{0}(\mathcal{A})$.
Hence $0\notin\sigma(M_wEM_u)=\sigma(EM_{uw})$.\\

(c) By Theorem 3.1 we have
$$\sigma_{p}(M_wEM_u)\setminus\{0\}=\sigma_{p}(EM_{u}M_{w})\setminus\{0\}=\sigma_{p}(EM_{uw})\setminus\{0\}.$$
So
$$\sigma_{p}(M_wEM_u)\setminus\{0\}=\{\lambda\in\mathbb{C}\setminus\{0\}:\mu(A_{\lambda,w})>0\}.$$\\

The proof of (d) and (e) are as the same as (c).

\vspace*{0.3cm} {\bf Corollary 3.4.} If
$M_wEM_u:L^2(\Sigma)\rightarrow L^2(\Sigma)$ and $|E(uw)|^2\geq
 E(|u|^2)E(|w|^2)$, then\\

 (a)
 $$\sigma_{jp}(M_wEM_u)\setminus \{0\}=ess \
 range(E(uw))\setminus\{0\}=\{\lambda\in\mathbb{C}\setminus\{0\}:\mu(A_{\lambda,w})>0\}.$$

(b) If $S\cap G=X$, then

 $$\sigma_{jp}(M_wEM_u)=ess \
 range(E(uw))=\{\lambda\in\mathbb{C}:\mu(A_{\lambda,w})>0\}.$$\\

 By Theorem 3.2 we have the following theorem.\\

\vspace*{0.3cm} {\bf Theorem 3.5.} For
$T=M_wEM_u:L^p(\Sigma)\rightarrow L^p(\Sigma)$ and $1\leq p\leq
\infty$ we have
$r(T)=\|E(uw)\|_{\infty}$.\\

Recall that an operator $T$ is quasinilpotent if
$\sigma(T)=\{0\}$. In light of Theorem 3.1, we have the following
corollary.\\

\vspace*{0.3cm} {\bf Corollary 3.6.} If
$T=M_wEM_u:L^p(\Sigma)\rightarrow L^p(\Sigma)$, $1\leq p\leq
\infty$ and $\frac{1}{p}+\frac{1}{p'}=1$, then the followings are
equivalent:\\

(a) $M_wEM_u$ is quasinilpotent;\\

(b) $E(uw)\equiv 0$;\\

(c) $M_{wE(uw)}EM_u\equiv 0$.\\

\vspace*{0.3cm} {\bf Proof.} $(a\Leftrightarrow b)$ Since
$\sigma(M_wEM_u)\setminus \{0\}=ess \ range (E(uw))\setminus
\{0\}$, it follows that $\sigma(M_wEM_u)=0$ if and only if
$E(uw)=0$ almost every where.\\

$(b\Leftrightarrow c)$ Since
$$\|M_{wE(uw)}EM_u\|=\|E(uw)(E(|w|^p))^{\frac{1}{p}}(E(|u|^{p'}))^{\frac{1}{p'}}\|_{\infty}$$
and $$|E(uw)|\leq
(E(|w|^p))^{\frac{1}{p}}(E(|u|^{p'}))^{\frac{1}{p'}},$$ it follows
that $E(uw)=0$, a.e, if and only if
$E(uw)(E(|w|^p))^{\frac{1}{p}}(E(|u|^{p'}))^{\frac{1}{p'}}=0$,
a.e, if and only if $\|M_{wE(uw)}EM_u\|=0$.

Now, we give some example of conditional expectation.\\

{\bf Example} (a) Let $X=\mathbb{N}\cup\{0\}$,
$\mathcal{G}=2^{\mathbb{N}}$ and let
$\mu(\{x\})=\frac{e^{-\theta}\theta^x}{x !}$, for each $x\in X$
and $\theta\geq0$. Elementary calculations show that $\mu$ is a
probability measure on $\mathcal{G}$. Let $\mathcal{A}$ be the
$\sigma$-algebra generated by the partition $B=\{\emptyset, X,
\{0\}, X_1=\{1, 3, 5, 7, 9, ....\}, X_2=\{2, 4, 6, 8, ....\},\}$
of $\mathbb{N}$. Note that, $\mathcal{A}$ is a
sub-$\sigma$-algebra of $\Sigma$ and each of element of
$\mathcal{A}$ is an $\mathcal{A}$-atom. Thus, the conditional
expectation of any $f\in \mathcal{D}(E)$ relative to $\mathcal{A}$
is constant on $\mathcal{A}$-atoms. Hence there exists scalars
$a_1, a_2, a_3$ such that

$$E(f)=a_1\chi_{0}+a_2\chi_{X_1}+a_3\chi_{X_2}.$$
So
$$E(f)(0)=a_1, \ \ \ \ E(f)(2n-1)=a_2, \ \ \ \ E(f)(2n)=a_3,$$
foe all $n\in \mathbb{N}$. By definition of conditional
expectation with respect to $\mathcal{A}$, we have
$$f(0)\mu(\{0\})=\int_{\{0\}}fd\mu=\int_{\{0\}}E(f)d\mu=a_1\mu(\{0\}),$$
so $a_1=f(0)$. Also,
$$\sum_{n\in \mathbb{N}}f(2n-1)\frac{e^{-\theta}\theta^{2n-1}}{(2n-1)
!}=\int_{X_1}fd\mu=\int_{X_1}E(f)d\mu$$$$=a_2\mu(X_2)=a_2\sum_{n\in
\mathbb{N}}\frac{e^{-\theta}\theta^{2n-1}}{(2n-1) !}.$$ and so

$$a_2=\frac{\sum_{n\in \mathbb{N}}f(2n-1)\frac{e^{-\theta}\theta^{2n-1}}{(2n-1)
!}}{\sum_{n\in \mathbb{N}}\frac{e^{-\theta}\theta^{2n-1}}{(2n-1)
!}}.$$ By the same method we have

$$a_3=\frac{\sum_{n\in \mathbb{N}}f(2n)\frac{e^{-\theta}\theta^{2n}}{(2n)
!}}{\sum_{n\in \mathbb{N}}\frac{e^{-\theta}\theta^{2n}}{(2n)
!}}.$$

If we set $f(x)=x$, then $E(f)$ is a special function as follows;

$$E(f)=\theta
coth(\theta)\chi_{X_1}+\frac{cosh(\theta)-1}{cosh(\theta)}\chi_{X_2}.$$\\

Also, if $u$ and $w$ are real functions on $X$ such that $M_wEM_u$
is boundede on $l^p$, then by Theorem 3.2 we have
$$\sigma(M_wEM_u)=\{u(0)w(0),\frac{\sum_{n\in \mathbb{N}}u(2n-1)w(2n-1)\frac{e^{-\theta}\theta^{2n-1}}{(2n-1)
!}}{\sum_{n\in \mathbb{N}}\frac{e^{-\theta}\theta^{2n-1}}{(2n-1)
!}},\frac{\sum_{n\in
\mathbb{N}}u(2n)w(2n)\frac{e^{-\theta}\theta^{2n}}{(2n)
!}}{\sum_{n\in \mathbb{N}}\frac{e^{-\theta}\theta^{2n}}{(2n)
!}}\}.$$

(b) Let $X=\mathbb{N}$, $\mathcal{G}=2^{\mathbb{N}}$ and let
$\mu(\{x\})=pq^{x-1}$, for each $x\in X$, $0\leq p\leq 1$ and
$q=1-p$. Elementary calculations show that $\mu$ is a probability
measure on $\mathcal{G}$. Let $\mathcal{A}$ be the
$\sigma$-algebra generated by the partition
$B=\{X_1=\{3n:n\geq1\}, X^{c}_1\}$ of $X$. So, for every $f\in
\mathcal{D}(E^{\mathcal{A}})$,

$$E(f)=\alpha_1\chi_{X_1}+\alpha_2\chi_{X^c_1}$$
and direct computations show that

$$\alpha_1=\frac{\sum_{n\geq1}f(3n)pq^{3n-1}}{\sum_{n\geq1}pq^{3n-1}}$$
and
$$\alpha_2=\frac{\sum_{n\geq1}f(n)pq^{n-1}-\sum_{n\geq1}f(3n)pq^{3n-1}}{\sum_{n\geq1}pq^{n-1}-\sum_{n\geq1}pq^{3n-1}}.$$

For example, if we set $f(x)=x$, then $E(f)$ is a special function
as follows;

$$\alpha_1=\frac{3}{1-q^3}, \ \ \ \ \
\alpha_2=\frac{1+q^6-3q^4+4q^3-3q^2}{(1-q^2)(1-q^3)}.$$\\
So, if $u$ and $w$ are real functions on $X$ such that $M_wEM_u$
is boundede on $l^p$, then by Theorem 3.2 we have
$$\sigma(M_wEM_u)=\{\frac{\sum_{n\geq1}u(3n)w(2n)pq^{3n-1}}{\sum_{n\geq1}pq^{3n-1}},
\frac{\sum_{n\geq1}u(n)w(n)pq^{n-1}-\sum_{n\geq1}u(3n)w(3n)pq^{3n-1}}{\sum_{n\geq1}pq^{n-1}-\sum_{n\geq1}pq^{3n-1}}\}.$$\\

(c) Let $X=[0,1)$, $\Sigma$=sigma algebra of Lebesgue measurable
subset of $X$, $\mu$=Lebesgue measure on $X$. Let $s:[0,
1)\rightarrow[0, 1)$ be defined by $s(x)=x+\frac{1}{4}$(mod 1).
Let $\mathcal{B}=\{E\in \Sigma:s(E)=E\}$. In this case
$$E^{\mathcal{B}}(f)(x)=\frac{f(x)+f(s(x))+f(s^{2}(x))+f(s^{3}(x))}{4},$$ where $s^j$
denotes the jth iteration of $s$. Also, $|f|\leq
3E^{\mathcal{B}}(|f|)$
a.e. Hence, the operator $EM_u$ is bounded on $L^p([0,1))$ if and only if $u\in L^{\infty}([0,1))$.\\
\\

(d) Let $X=[0,a]\times [0,a]$ for $a>0$, $d\mu=dxdy$, $\Sigma$ the
Lebesgue subsets of $X$ and let $\mathcal{A}=\{A\times [0,a]: A \
\mbox{is a Lebesgue set in} \ [0,a]\}$. Then, for each $f\in
\mathcal{D}(E)$, $(Ef)(x, y)=\int_0^af(x,t)dt$, which is
independent of the second coordinate. For example, if we set
$a=1$, $w(x,y)=1$ and $u(x,y)=e^{(x+y)}$, then
$E(u)(x,y)=e^x-e^{x+1}$ and $M_wEM_u$ is bounded. Therefore by
Theorem 3.2 $\sigma(M_wEM_u)=[e-e^2, 1-e]$.

%Now we
%draw the graph of $f$ and $E(f)$ for $a=1$, so one can see the
%effect of
%$E$ on $f$.\\
%
%\includegraphics[width=9.5 cm, height=8cm]{estaremi5.pdf}
%
%\includegraphics[width=9.5 cm, height=8cm]{estaremi4.pdf}\\


\begin{thebibliography}{99}


\bibitem{dhd} P.G. Dodds, C.B. Huijsmans and B. De Pagter,
characterizations of conditional expectation-type operators,
Pacific J. Math. {\bf 141}(1) (1990), 55-77.


\bibitem{e} Y. Estaremi, Essential norm of weighted conditional type operators on $L^p$-spaces, to appear in positivity.

\bibitem{ej1} Y. Estaremi and M.R. Jabbarzadeh, Weighted lambert type operators on
$L^{p}$-spaces, Oper. Matrices {\bf 1} (2013), 101-116.

\bibitem{ej2} Y. Estaremi and M.R. Jabbarzadeh, Weighted composition lambert type operators on
$L^{p}$-spaces, to appear in Mediterranean Journal of Mathematics.



\bibitem{g}
J. J. Grobler and B. de Pagter, Operators representable as
multiplication-conditional expectation operators, J. Operator
Theory {\bf 48} (2002), 15-40.


\bibitem{her}
J. Herron, Weighted conditional expectation operators, Oper.
Matrices {\bf 1} (2011), 107-118.

\bibitem{lam}
A. Lambert, $L^p$ multipliers and nested sigma-algebras, Oper.
Theory Adv. Appl. {\bf 104} (1998),  147-153.


\bibitem{mo}
Shu-Teh Chen, Moy,  Characterizations of conditional expectation
as a transformation on function spaces, Pacific J. Math. {\bf 4}
(1954), 47-63.
\bibitem{rao}
M. M. Rao, Conditional measure and applications, Marcel Dekker,
New York, 1993.

\bibitem{ty}
T. Yamazaki, An expression of spectral radius via Aluthge
transformation, Pro. A. M. S. {\bf 4} (2001), 1131-1137.





\end{thebibliography}
\end{document}